\documentclass[12pt]{article}
\usepackage{amssymb}
\usepackage{mathrsfs}
\makeatletter
 
 \@addtoreset{equation}{section}
\makeatother

\newcommand{\rem}{\noindent{\bf Remark. }}

\newcommand{\qed}{\hbox{\rule[-2pt]{5pt}{11pt}}}

\newtheorem{dfn}{Definition}[section]
\newtheorem{thm}[dfn]{Theorem}
\newtheorem{prop}[dfn]{Proposition}

\newtheorem{conj}[dfn]{Conjecture}

\begin{document}
%\begin{linenumbers}
%\begin{frontmatter}
\title{On some families of invariant polynomials divisible by three and their zeta functions}
\author{Koji Chinen\footnotemark[1]}
\date{}
\maketitle

\begin{abstract}
In this note, we establish an analog of the Mallows-Sloane bound for Type III formal weight enumerators. This completes the bounds for all types (Types I through IV) in synthesis of our previous results. Next we show by using the binomial moments that there exists a family of polynomials divisible by three, which are not related to linear codes but are invariant under the MacWilliams transform for the value 3/2. We also discuss some properties of the zeta functions for such polynomials. 
\end{abstract}%  words

\footnotetext[1]{Department of Mathematics, School of Science and Engineering, Kindai University. 3-4-1, Kowakae, Higashi-Osaka, 577-8502 Japan. E-mail: chinen@math.kindai.ac.jp}
\footnotetext{This work was supported by JSPS KAKENHI Grant Number JP26400028. }

%\footnotetext[1]{}

%%%%%%%%%%%%%%%%%%%%%%%%%%%%%%%%%%%%%%%%%%%%%%%%%%%%%%%%%%%%%
%%%%%%%%%%%%%%%%%%%%%% section 1 %%%%%%%%%%%%%%%%%%%%%%%%%%%%
\section{Introduction}\label{section:intro}
This article, as a sequel of \cite{Ch3}-\cite{Ch5}, investigates some invariant polynomials which are divisible by three, that is, polynomials of the form
\begin{equation}\label{eq:fwe_form} %%%%%%%%%%%%%%%%%%% fwe_form
W(x,y)=x^n+\sum_{i=d}^n A_i x^{n-i} y^i\in {\bf C}[x,y]\quad (A_d\ne 0)
\end{equation}
that satisfy ``$A_i\ne 0$ $\Rightarrow$ $3|i$'' and certain transformation rules: for a linear transformation $\sigma=\left(\begin{array}{cc} a & b \\ c & d \end{array}\right)$, the action of $\sigma$ on $W(x,y)$ is defined by
$$W^\sigma(x,y)=W(ax+by, cx+dy)$$
and we are interested in $W(x,y)$ of the form (\ref{eq:fwe_form}) with the property 
$$W^{\sigma_q}(x,y)=\pm W(x,y),$$
where 
$$\sigma_q=\frac{1}{\sqrt{q}}\left(\begin{array}{rr} 1 & q-1 \\ 1 & -1 \end{array}\right) \quad(\mbox{the MacWilliams transform}).$$
We call $W(x,y)$ with $W^{\sigma_q}(x,y)=W(x,y)$ a ``$\sigma_q$-invariant polynomial'' and $W(x,y)$ with $W^{\sigma_q}(x,y)=-W(x,y)$ a ``formal weight enumerator for $q$''. In general, $W(x,y)$ is called ``divisible by $c$'' ($c>1$) if ``$A_i\ne 0$ $\Rightarrow$ $c|i$''. In what follows, we put $\tau=\left(\begin{array}{cc} 1 & 0 \\ 0 & \omega \end{array}\right)$ ($\omega=(-1+\sqrt{-3})/2$). The Pochhammer symbol $(a)_n$ means $(a)_n=a(a+1)\cdots (a+n-1)$ for $n\geq 1$ and $(a)_0=1$. 

The earliest example of the divisible formal weight enumerator in the literature is the case $(q,c)=(2,4)$, which is given in Ozeki \cite{Oz} (the denomination ``formal weight enumerator'' is also due to him). Ozeki's formal weight enumerators are members of the polynomial ring
$$R_{\rm II}^- :={\bf C}[W_{{\mathcal H}_{8}}(x,y), W_{12}(x,y)],$$
where
\begin{eqnarray*}
W_{{\mathcal H}_{8}}(x,y) &=& x^8+14x^4y^4+y^8, \\%\label{eq:we_hamming}\\
W_{12}(x,y) &=& x^{12}-33x^8y^4-33x^4y^8+y^{12} %\label{eq:fwe_deg12}
\end{eqnarray*}
($W_{12}(x,y)$ satisfies ${W_{12}}^{\sigma_2}(x,y)=-W_{12}(x,y)$). Note that $W_{{\mathcal H}_{8}}(x,y)$ is the weight enumerator of the extended Hamming code. We will call formal weight enumerators in $R_{\rm II}^-$ ``Type II formal weight enumerators'', since they resemble weight enumerators of Type II codes, which are divisible by four and $\sigma_2$-invariant. We also have rings of formal weight enumerators for the cases $(q,c)=(2,2), (3,3)$ and $(4,2)$ which we shall call Types I, III and IV, respectively: 
\begin{eqnarray*}
(\mbox{Type I}) & & R_{\rm I}^- :={\bf C}[W_{2,2}(x,y), \varphi_4(x,y)],\\
(\mbox{Type III}) & & R_{\rm III}^- :={\bf C}[W_4(x,y), \psi_6(x,y)],\\
(\mbox{Type IV}) & & R_{\rm IV}^- :={\bf C}[W_{2,4}(x,y), \varphi_3(x,y)],
\end{eqnarray*}
where, 
\begin{eqnarray*}
W_{2,q}(x,y) &=& x^2+(q-1)y^2, \\
\varphi_4(x,y) &=& x^4-6x^2y^2+y^4, \\
W_4(x,y) &=& x^4+8xy^3, \\
\psi_6(x,y) &=& x^6-20x^3y^3-8y^6, \\
\varphi_3(x,y) &=& x^3-9xy^2.
\end{eqnarray*}
The ring $R_{\rm III}^-$ is introduced by Ozeki \cite{Oz2}, $R_{\rm I}^-$ and $R_{\rm IV}^-$ are dealt with in \cite{Ch5}. 

Our first goal in this article is to complete the following theorem by proving the case of Type III (the cases Types I and IV are already proved in \cite{Ch5} and the case Type II is proved in \cite{Ch1}): 
%%%%%%%%%%%%%%%%% Theorem (analog of Mallows-Sloane) %%%%%%%%%%%%%%%%%%%
\begin{thm}\label{thm:analog_MS}
For all formal weight enumerators of Types I through IV of the form (\ref{eq:fwe_form}), we have the following:
\begin{eqnarray}
(\mbox{\rm Type I}) & & d\leq 2\left[ \frac{n-4}{8} \right] +2,\nonumber\\
(\mbox{\rm Type II}) & & d\leq 4\left[ \frac{n-12}{24} \right] +4,\nonumber\\
(\mbox{\rm Type III}) & & d\leq 3\left[ \frac{n-6}{12} \right] +3,\label{eq:analogMS_III}\\
(\mbox{\rm Type IV}) & & d\leq 2\left[ \frac{n-3}{6} \right] +2,\nonumber
\end{eqnarray}
where $[x]$ means the greatest integer not exceeding $x$ for $x\in{\bf R}$. 
\end{thm}
%%%%%%%%%%%%%%%%%%%%%%%%%%%%%%%%%%%%%%%%%%%%%%%%
This is an analog of the famous Mallows-Sloane bound for weight enumerators of divisible self-dual codes (\cite{MalSl}). Similarly to the case of codes, we can define the extremal formal weight enumerator: 
%%%%%%%%%%%%%%%%% Definition (extremal FWE) %%%%%%%%%%%%%%%%%%%
\begin{dfn}\label{dfn:extremalFWE}
A formal weight enumerator of Types I through IV is called extremal if an equality holds in Theorem \ref{thm:analog_MS}. 
\end{dfn}
%%%%%%%%%%%%%%%%%%%%%%%%%%%%%%%%%%%%%%%%%%%%%%%%
Our interest in divisible formal weight enumerators arose from the consideration of their zeta functions. Zeta functions of this type were defined in Duursma \cite{Du1} for weight enumerators of linear codes (see also \cite{Du2}-\cite{Du4}) and some generalization was made by the present author (\cite{Ch1}, \cite{Ch2}): 
%%%%%%%%%%%%%%%%% Definition (zeta function) %%%%%%%%%%%%%%%%%%%
\begin{dfn}\label{dfn:zeta}
For any homogeneous polynomial of the form (\ref{eq:fwe_form}) and $q\in{\bf R}$ ($q>0, q\ne 1$), there exists a unique polynomial $P(T)\in{\bf C}[T]$ of degree at most $n-d$ such that
\begin{equation}\label{eq:zeta_duursma}
\frac{P(T)}{(1-T)(1-qT)}(y(1-T)+xT)^n=\cdots +\frac{W(x,y)-x^n}{q-1}T^{n-d}+ \cdots.
\end{equation}
We call $P(T)$ and $Z(T)=P(T)/(1-T)(1-qT)$ the zeta polynomial and the zeta function of $W(x,y)$, respectively. 
\end{dfn}
%%%%%%%%%%%%%%%%%%%%%%%%%%%%%%%%%%%%%%%%%%%%%%%%
We must assume $d, d^\perp \geq 2$ where $d^\perp$ is defined by 
$$W^{\sigma_q}(x,y)=\pm x^n + A_{d^\perp} x^{n-d^\perp} y^{d^\perp}+ \cdots,$$
when considering zeta functions (\cite[p.57]{Du2}). The Riemann hypothesis is formulated as follows: 
%%%%%%%%%%%%%%%%%% Definition (RH) %%%%%%%%%%%%%%%%%%
\begin{dfn}[Riemann hypothesis]\label{dfn:RH}
A polynomial of the form (\ref{eq:fwe_form}) with $W^{\sigma_q}(x,y)=\pm W(x,y)$ satisfies the Riemann hypothesis if all the zeros of $P(T)$ have the same absolute value $1/\sqrt{q}$. 
\end{dfn}
Our second result is the following theorem, which is an analog of Okuda's theorem (\cite[Theorem 5.1]{Ok}), of which proof will be given briefly in Section \ref{section:typeIII} : 
%%%%%%%%%%%%%%%%% Theorem (analog Okuda) %%%%%%%%%%%%%%%%%%%
\begin{thm}\label{thm:analog_okuda}
Let $W(x,y)$ be the Type III extremal formal weight enumerator of degree $n=12k+6$ ($k\geq 1$). Then 
$$W^\ast (x,y):=\frac{1}{(n-3)_4}\frac{\partial}{\partial x}\left( \frac{\partial^3}{\partial x^3} +  \frac{\partial^3}{\partial y^3}\right)W(x,y)$$
is the extremal formal weight enumerator of degree $n-4$. Moreover, the zeta polynomial $P(T)$ of $W(x,y)$ and $P^\ast(T)$, that of $W^\ast(x,y)$ are related by $P^\ast(T)=(3T^2-3T+1)P(T)$. The Riemann hypothesis of $W(x,y)$ and that of $W^\ast(x,y)$ are equivalent. 
\end{thm}
%%%%%%%%%%%%%%%%%%%%%%%%%%%%%%%%%%%%%%%%%%%%%%%%
These results, together with the ones in \cite{Ch1} and \cite{Ch5} suggest that formal weight enumerators of Types I through IV have similar properties to the weight enumerators of corresponding Types. 

The last feature of this article is the discovery of $\sigma_{3/2}$-invariant polynomials. They are also divisible by three: 
\begin{equation}\label{eq:ring_inv3over2}
R_{3/2}:={\bf C}[\eta_6(x,y), \eta_{24}(x,y)],
\end{equation}
where
\begin{eqnarray}
\eta_6(x,y) &=& x^6+\frac{5}{2}x^3y^3-\frac{1}{8}y^6, \label{eq:3over2_deg6}\\
\eta_{24}(x,y) &=& x^{24}+\frac{253}{4}x^{18}y^6+\frac{1265}{32}x^{15}y^9+\frac{7659}{256}x^{12}y^{12}\nonumber\\
          & & -\frac{1265}{256}x^9y^{15}+\frac{253}{256}x^6y^{18}+\frac{1}{4096}y^{24}. \label{eq:3over2_deg24}
\end{eqnarray}
We can also construct the ring of formal weight enumerators for $(q,c)=(3/2, 3)$: 
\begin{equation}\label{eq:ring_fwe3over2}
R_{3/2}^-:={\bf C}[\eta_6(x,y), \eta_{12}(x,y)],
\end{equation}
where
\begin{equation}\label{eq:3over2fwedeg6}
\eta_{12}(x,y) = x^{12}-11x^9y^3-\frac{11}{8}x^3y^9-\frac{1}{64}y^{12}.
\end{equation}
These families were discovered by the use of the binomial moments. We will explain it and observe their Riemann hypothesis in Section \ref{section:3over2}. 
%%%%%%%%%%%%%%%%%%%%%%%%%%%%%%%%%%%%%%%%%%%%%%%%%%%%%%%%%%%%%
%%%%%%%%%%%%%%%%%%%%%% section 2 %%%%%%%%%%%%%%%%%%%%%%%%%%%%
\section{Type III formal weight enumerators}\label{section:typeIII}
First we give an outline of the proof of Theorem \ref{thm:analog_MS} (Type III). For a homogeneous polynomial $p(x,y)\in {\bf C}[x,y]$, $p(x,y)(D)$ means a differential operator obtained by replacing $x$ by $\partial/\partial x$ and $y$ by $\partial/\partial y$. Here we use $p(x,y)=y(y^3-8x^3)$. In a similar manner to Duursma \cite[Lemma 2]{Du4}, we can prove the following (see also \cite[Proposition 3.1]{Ch5}): 
%%%%%%%%%%%%%%%%%% proposition %%%%%%%%%%%%%%%%%%
\begin{prop}\label{prop:dividsion}
Let $W(x,y)$ be a Type III formal weight enumerator with $d\geq 6$. Then we have 
\begin{equation}\label{eq:divide_III}
\{y(x^3-y^3)\}^{d-4}| p(x,y)(D)W(x,y).
\end{equation}
\end{prop}
Using this, we can prove (\ref{eq:analogMS_III}). Proof is similar to that of \cite[Theorem 3.3]{Ch5} and the notation follows it: 

\medskip
\noindent{\bf (Proof of Theorem \ref{thm:analog_MS} (Type III))}

\medskip
Let $a(x,y)=\{y(x^3-y^3)\}^{d-4}$ and we put 
$$p(x,y)(D)W(x,y)=a(x,y)\tilde a(x,y).$$
Then from the transformation rules of $p(x,y)$, $a(x,y)$ and $W(x,y)$ by $\sigma_3$, $^t\!\sigma_3$ and $\tau$, we can show that $\tilde a(x,y)$ is a constant times a formal weight enumerator (the discussion is similar to that of \cite[Theorem 3.3]{Ch5}). So we have $\psi_6(x,y)|\tilde a(x,y)$. Since $(a(x,y), \psi_6(x,y))=1$, we can conclude 
$$a(x,y)\psi_6(x,y)|p(x,y)(D)W(x,y).$$
Comparing the degrees on both sides, we have $4(d-4)+6 \leq n-4$. Putting $d=3d'$ ($d\in{\bf N}$), we get $d'\leq (n-6)/12+1$. We can obtain the conclusion immediately. \qed

\medskip
\noindent\rem Some numerical examples of zeta polynomials for Type III formal weight enumerators are given and the extremal property is mentioned up to degree 18 in \cite[Section 4]{Ch1}. 

\medskip
\noindent{\bf (Proof of Theorem \ref{thm:analog_okuda})}

\medskip
We follow the method of Okuda \cite{Ok} (see also \cite[Theorem 3.9]{Ch5}). Here we use $p(x,y)=x(x^3+y^3)$. Let $W(x,y)$ be the extremal formal weight enumerator of degree $n=12k+6$. Then from the rules 
\begin{eqnarray*}
& & p^{^t \sigma_3}(x,y)=p^{^t \tau}(x,y)=p(x,y),\\
& & W^{\sigma_3}(x,y)=-W(x,y),\ W^{\tau}(x,y)=W(x,y)
\end{eqnarray*}
and the uniqueness of the extremal formal weight enumerator at each degree, we can see that $W^\ast(x,y)$ is the extremal formal weight enumerator of degree $12k+2$ (use \cite[Theorem 2.3 (i)]{Ch5}). Next, by the use of the MDS weight enumerators for $q=3$, we can see that the zeta polynomial of $x^4(D) W(x,y)/(n-3)_4$ is $P(T)$, that of $xy^3(D) W(x,y)/(n-3)_4$ is $(1-T)^3P(T)$. Adjusting the degrees, we can conclude that $P^\ast(T)=(3T^2-3T+1)P(T)$. The equivalence of the Riemann hypothesis is straightforward. \qed
%%%%%%%%%%%%%%%%%%%%%%%%%%%%%%%%%%%%%%%%%%%%%%%%%%%%%%%%%%%%%
%%%%%%%%%%%%%%%%%%%%%% section 3 %%%%%%%%%%%%%%%%%%%%%%%%%%%%
\section{Polynomials for $q=3/2$}\label{section:3over2}
Our construction of $\eta_6(x,y)$ (see (\ref{eq:3over2_deg6})) uses the binomial moments. The method is much the same as that of \cite{Ch3}, so we give an outline. We search a $\sigma_q$-invariant polynomial $W(x,y)$ of the form
\begin{equation}\label{eq:div3general}
W(x,y)=\sum_{i=0}^{[2n/3]}A_i x^{2n-3i}y^{3i}\quad (A_0=1).
\end{equation}
The formula of the binomial moments for (\ref{eq:div3general}) becomes
\begin{equation}\label{eq:moment_div3}
\sum_{i=0}^{[(2n-\nu)/3]} {{2n-3i}\choose{\nu}} A_i - q^{n-\nu} \sum_{i=0}^{[\nu/3]} {{2n-3i}\choose{2n-\nu}} A_i =0
\quad (\nu=0, 1, \cdots, 2n)
\end{equation}
(it is easily obtained from \cite[p.131, Problem (6)]{MaSl}). In (\ref{eq:moment_div3}), the values $\nu$ and $2n-\nu$ give essentially the same formula, so it suffices to consider the cases $\nu=0, 1, \cdots, n$. Moreover, (\ref{eq:moment_div3}) is trivial when $\nu=n$. Thus (\ref{eq:moment_div3}) gives $n$ linear equations of $[2n/3]+1$ unknowns $A_0, A_1, \cdots, A_{[2n/3]}$. The number of equations and unknowns coincide when $n=3$, in which case the system of equations becomes
$$\left\{\begin{array}{l}
(1-q^3)A_0+A_1+A_2=0,\\
6(1-q^2)A_0+3A_1=0,\\
15(1-q)A_0+3A_1=0.
\end{array}\right.$$
Since $A_0=1$, we have $2q^2-5q+3=0$. We get a non-trivial value $q=3/2$. We can determine other coefficients $A_1=5/2$, $A_2=-1/8$ and get $\eta_6(x,y)$. We can verify it is indeed $\sigma_{3/2}$-invariant. We can also verify (with some computer algebra system) that there is no $\sigma_{3/2}$-invariant polynomial $W(x,y)$ of even degrees in the range $8\leq \deg W(x,y) \leq 22$ except for $\eta_6(x,y)^2$ and $\eta_6(x,y)^3$, but we can find $\eta_{24}(x,y)$ in (\ref{eq:3over2_deg24}) at degree 24 ($\eta_6(x,y)$ and $\eta_{24}(x,y)$ are algebraically independent). We can furthermore find $\eta_{12}(x,y)$ from the condition that it is invariant under $\sigma_{3/2}\tau\sigma_{3/2}$. The ring $R_{3/2}$ is the invariant polynomial ring of the group $\langle \sigma_{3/2}, \tau \rangle$, $R_{3/2}^-$ is that of $\langle \sigma_{3/2}\tau\sigma_{3/2}, \tau \rangle$. 

For the members of $R_{3/2}^-$ (including $R_{3/2}$), there seems to be bounds similar to Theorem \ref{thm:analog_MS} (proof seems to be difficult): 
%%%%%%%%%%%%%%%%% Theorem (analog of Mallows-Sloane 3/2) %%%%%%%%%%%%%%%%%%%
\begin{conj}\label{conj:analog_MS_3over2}
(i) All $\sigma_{3/2}$-invariant polynomials of the form (\ref{eq:fwe_form}) in $R_{3/2}$ satisfy
$$d\leq 3\left[ \frac{n}{24} \right] +3.$$
(ii) All formal weight enumerators of the form (\ref{eq:fwe_form}) in $R_{3/2}^-$ satisfy
$$d\leq 3\left[ \frac{n-12}{24} \right] +3.$$
\end{conj}
%%%%%%%%%%%%%%%%%%%%%%%%%%%%%%%%%%%%%%%%%%%%%%%%
Here are some examples of zeta polynomials for the members of $R_{3/2}^-$. The zeta polynomial of $\eta_6(x,y)$ is $P_6(T)=(3T^2+3T+2)/8$, that of $\eta_{12}(x,y)$ is $P_{12}(T)=(3T^2-2)(27T^6+27T^5+36T^4+26T^3+24T^2+12T+8)/160$ (the zeta polynomial of $\eta_{24}(x,y)$ is a polynomial of degree 14). From numerical experiments we can conjecture that extremal $\sigma_{3/2}$-invariant polynomials and extremal formal weight enumerators in $R_{3/2}^-$ satisfy the Riemann hypothesis. 

%\end{linenumbers}
\end{document}